\documentclass[12pt]{article}
\usepackage{amssymb}

\newtheorem{theorem}{{\bf Theorem}}[section]
\newtheorem{lemma}{{\bf Lemma}}[section]
\newtheorem{proposition}{{\bf Proposition}}[section]
\newtheorem{definition}{{\bf Definition}}[section]

\newtheorem{remark}{{\bf Remark}}[section]

\input{tcilatex}
\input tcilatex
\begin{document}

\title{Reflected BSDE with monotonicity and general increasing in $y$, and non-Lipschitz
conditions in $z$}
\author{Mingyu Xu \thanks{Email: xvmingyu@gmail.com}\\
{\small Departement des Math\'ematiques, Universit\'e du Maine,
France;}\\{\small Department of Financial Mathematics and Control
science, School
of Mathematical Science,} \\
{\small Fudan University, Shanghai, 200433, China.}\\
Draft : March 1st, 2006} \maketitle

\textbf{Abstract} In this paper, we study the reflected BSDE with one
continuous barrier, under the monotonicity and general increasing condition
on $y$ and non Lipschitz condition on $z$. We prove the existence and
uniqueness of the solution to these equation by approximation method.


\section{Introduction}

Nonlinear backward stochastic differential equations (BSDE in short) were
firstly introduced by Pardoux and Peng in 1990, \cite{PP90}. They proved
that there exists a unique solution $(Y,Z)$ to this equation if the terminal
condition $\xi $ and coefficient $f$ satisfy smooth square-integrability
assumptions and $f(t,\omega ,y,z)$ is Lipschitz in $(y,z)$ uniformy in $%
(t,\omega )$. Later many assumptions have been made to relax the Lipschitz
condition on $f$. Pardoux (1999, \cite{P99}) and Briand et al. (2003, \cite
{BDHPS}) studied the solution of a BSDE with a coefficient $f(t,\omega ,y,z)$%
, which still satisfies the Lipschitz condition on $z$, but only
monotonicity, continuity and generalized increasing on $y$, i.e.for some
continuous increasing function $\varphi :\mathbb{R}_{+}\rightarrow \mathbb{R}%
_{+}$, real number $\mu $ $>0$:
\begin{eqnarray}
\left| f(t,y,0)\right| &\leq &\left| f(t,0,0)\right| +\varphi
(\left| y\right| )\mbox{, }\forall (t,y)\in [0,T]\times
\mathbb{R}\mbox{, a.s.;}
\label{con-mono} \\
(y-y^{\prime })(f(t,y,z) &-&f(t,y^{\prime },z))\leq \mu (y-y^{\prime })^{2}%
\mbox{, }\forall (t,z)\in [0,T]\times \mathbb{R}^{d}\mbox{,
}y,y^{\prime }\in \mathbb{R}\mbox{, a.s.}  \nonumber
\end{eqnarray}
The case when $f$ is quadratic on $z$ and $\xi $ is bounded was firstly
studied by Kobylanski in \cite{K00}. She proved an existence result when the
coefficient is only linear growth in $y$, and quadratic in $z$. In \cite
{LS98}, Lepeltier and San Mart\'{i}n generalized to a superlinear case in $y$%
. More recently, in \cite{LS04}, they and Briand considered the BSDE whose
coefficient $f$ satisfies only monotonicity, continuity and generalized
increasing on $y$, and quadratic or linear increasing in $z$, i.e.
\begin{eqnarray}
(y-y^{\prime })(f(t,y,z) &-&f(t,y^{\prime },z))\leq \mu (y-y^{\prime })^{2}%
\mbox{, }\forall (t,z)\in [0,T]\times \mathbb{R}^{d}\mbox{,
}y,y^{\prime
}\in \mathbb{R}\mbox{, a.s.}  \nonumber \\
\left| f(t,y,z)\right| &\leq &\varphi (\left| y\right| )+A\left|
z\right| ^{2}\mbox{, }\forall (t,y)\in [0,T]\times \mathbb{R}\mbox{,
a.s.;} \label{con-qua}
\end{eqnarray}
or
\begin{equation}
\left| f(t,y,z)\right| \leq g_{t}+\varphi (\left| y\right| )+A\left|
z\right| ,\forall (t,y)\in [0,T]\times \mathbb{R}\mbox{, a.s..}
\label{con-lin}
\end{equation}
In the same paper, they studied the case $f(t,y,z)=\left| z\right| ^{p}$,
for $p\in (1,2]$, and gave some sufficient and necessary conditions on $\xi $
for the existence of solutions.

El Karoui, Kapoudjian, Pardoux, Peng and Quenez introduced the notion of
reflected BSDE (RBSDE in short) on one lower barrier in 1997, \cite{EKPPQ}:
the solution is forced to remain above a continuous process, which is
considered as the lower barrier. More precisely, a solution for such
equation associated to a coefficient $f(t,\omega ,y,z)$, a terminal value $%
\xi $, a continuous barrier $L$, is a triple
$(Y_{t},Z_{t},K_{t})_{0\leq t\leq T}$ of adapted processes valued
on $\mathbb{R}^{1+d+1}$, which satisfies a square integrability
condition,
\[
Y_{t}=\xi
+\int_{t}^{T}f(s,Y_{s},Z_{s})ds+K_{T}-K_{t}-\int_{t}^{T}Z_{s}dB_{s},0\leq
t\leq T\mbox{, a.s.,}
\]
and $Y_{t}\geq L_{t}$, $0\leq t\leq T$, a.s.. Furthermore, the process $%
(K_{t})_{0\leq t\leq T}$ is non decreasing, continuous, and the role of $%
K_{t}$ is to push upward the state process in a minimal way, to keep it
above $L$. In this sense it satisfies $\int_{0}^{T}(Y_{s}-L_{s})dK_{s}=0$.
They proved the existence and uniqueness of the solution when $f(t,\omega
,y,z)$ is Lipschitz in $(y,z)$ uniformly in $(t,\omega )$. Then Matoussi
(1997, \cite{M97}) consider RBSDE's where the coefficient $f$ is continuous
and at most linear growth in $y$, $z$. In this case, he proved the existence
of maximal solution for the RBSDE.

In \cite{KLQT}, Kobylanski, Lepeltier, Quenez and Torres proved the
existence of a maximal and minimal bounded solution for the RBSDE when the
coefficient $f(t,\omega ,y,z)$ is super linear increasing in $y$ and
quadratic in $z$, i.e. there exists a function $l$ strictly positive such
that
\[
\left| f(t,y,z)\right| \leq l(y)+A\left| z\right| ^{2}\mbox{, with }%
\int_{0}^{\infty }\frac{dx}{l(x)}=+\infty .
\]
In this case, $\xi $ and $L$ are required to be bounded, and $L$ is a
continuous process. Recently, in \cite{LMX} Lepeltier, Matoussi and Xu
considered the case when $f(t,\omega ,y,z)$ satisfies (\ref{con-mono}) and
is Lipschitz in $z$. They proved the existence and uniqueness of the
solution by an approximation procedure.

In this paper, we study the RBSDEs whose the coefficient $f\ $satisfies the
conditions (\ref{con-qua}) or (\ref{con-lin}), when the lower barrier $L$ is
uniformly bounded. We prove the existence of a solution, following the
methods in \cite{LS04}, and we give a necessary and sufficient condition for
the case when $f(t,\omega ,y,z)=\left| z\right| ^{2}$, and its explicit
solution.

The paper is organized as follows: in Section 2, we present the basic
assumptions and the definition of the RBSDE; then in Section 3, we prove the
existence of a solution when $f(t,\omega ,y,z)$ satisfies the conditions (%
\ref{con-qua}), $\xi $ and $L$ are bounded; in the following section, we
consider the case when $f(t,\omega ,y,z)=\left| z\right| ^{2}$, and $\xi $
is not necessarily bounded. In this section, we give a necessary and
sufficient condition on the terminal condition $\xi $ for $p=2$ and its
explicit solution. Finally, in section 5, we study the RBSDE with the
condition (\ref{con-lin}), and prove the existence of a solution. At last,
in Appendix, we generalize the comparison theorem in \cite{KLQT}, and get
some comparison theorems, which help us to pass to the limit in the
approximations.

\section{Notations}

Let $(\Omega ,\mathcal{F},P)$ be a complete probability space, and $%
(B_{t})_{0\leq t\leq T}=(B_{t}^{1},B_{t}^{2},\cdots ,B_{t}^{d})_{0\leq t\leq
T}^{\prime }$ be a $d$-dimensional Brownian motion defined on a finite
interval $[0,T]$, $0<T<+\infty $. Denote by $\{\mathcal{F}_{t};0\leq t\leq
T\}$ the standard filtration generated by the Brownian motion $B$, i.e. $%
\mathcal{F}_{t}$ is the completion of
\[
\mathcal{F}_{t}=\sigma \{B_{s};0\leq s\leq t\},
\]
with respect to $(\mathcal{F},P)$. We denote by $\mathcal{P}$ the $\sigma $%
-algebra of predictable sets on $[0,T]\times \Omega $.

We will need the following spaces:

\[
\begin{array}{ll}
\mathbf{L}^{2}(\mathcal{F}_{t})= & \{\eta
:\mathcal{F}_{t}\mbox{-measurable
random real-valued variable, s.t. }E(|\eta |^{2})<+\infty \}, \\
\mathbf{H}_{n}^{2}(0,T)= & \{(\psi _{t})_{0\leq t\leq
T}:\mbox{predictable process valued in }\mathbb{R}^{n}\mbox{, s.t.
}E\int_{0}^{T}\left| \psi
(t)\right| ^{2}dt<+\infty \}, \\
\mathbf{S}^{2}(0,T)= & \{(\psi _{t})_{0\leq t\leq
T}:\mbox{progressively
measurable, continuous, real-valued process,} \\
& \mbox{s.t. }E(\sup_{0\leq t\leq T}\left| \psi (t)\right|
^{2})<+\infty \},
\\
\mathbf{A}^{2}(0,T)= & \{(K_{t})_{0\leq t\leq T}:\ \mbox{adapted
continuous
increasing process, } \\
& \mbox{s.t. }K(0)=0\mbox{, }E(K(T)^{2})<+\infty \}.
\end{array}
\]

Now we introduce the definition of the solution of reflected backward
stochastic differential equation with a terminal condition $\xi $, a
coefficient $f$ and a continuous reflecting lower barrier $L$(in short RBSDE$%
(\xi ,f,L)$), which is the same as in El Karoui et al.(1997, \cite{EKPPQ}).

\begin{definition}
\label{Def}We say that the triple $(Y_{t},Z_{t},K_{t})_{0\leq t\leq T}$ of
progressively measurable processes is a solution of RBSDE$(\xi ,f,L)$, if
the followings hold:

(i) $(Y_{t})_{0\leq t\leq T}\in \mathbf{S}^{2}(0,T)$, $(Z_{t})_{0\leq t\leq
T}\in \mathbf{H}_{d}^{2}(0,T)$, and $(K_{t})_{0\leq t\leq T}\in \mathbf{A}%
^{2}(0,T)$.

(ii) $Y_{t}=\xi
+\int_{t}^{T}f(s,Y_{s},Z_{s})ds+K_{T}-K_{t}-\int_{t}^{T}Z_{s}dB_{s},\;\;0%
\leq t\leq T$ a.s.

(iii) $Y_{t}\geq L_{t},\;\;0\leq t\leq T.$

(iv) $\int_{0}^{T}(Y_{s}-L_{s})dK_{s}=0,$ a.s.
\end{definition}

\section{The general case of $f$ quadratic increasing}

In this section, we work under the following assumptions:

\textbf{Assumption 1.} $\xi $ is an $\mathcal{F}_{T}$-adapted and bounded
random variable;

\textbf{Assumption 2.} a coefficient $f:\Omega \times [0,T]\times \mathbb{%
R\times R}^{d}\rightarrow \mathbb{R}$, is such that for some
continuous increasing function $\varphi :\mathbb{R}_{+}\rightarrow
\mathbb{R}_{+}$, real numbers $\mu $ and $A>0$ and $\forall
(t,y,y^{\prime }z)\in [0,T]\times \mathbb{R\times R\times R}^{d}$,

$
\begin{array}{cl}
\mbox{(i)} & f(\cdot ,y,z)\mbox{ is progressively measurable;} \\
\mbox{(ii)} & \left| f(t,y,z)\right| \leq \varphi (\left| y\right|
)+A\left|
z\right| ^{2}\mbox{;} \\
\mbox{(iii)} & (y-y^{\prime })(f(t,y,z)-f(t,y^{\prime },z))\leq \mu
(y-y^{\prime })^{2}\mbox{;} \\
\mbox{(iv)} & y\rightarrow f(t,y,z)\mbox{ is continuous, a.s.}
\end{array}
$

\textbf{Assumption 3.} a barrier $(L_{t})_{0\leq t\leq T}$, is a bounded
continuous progressively measurable real-valued process, $b:=\sup_{0\leq
t\leq T}\left| L_{t}\right| <+\infty $, $L_{T}\leq \xi $, a.s.

Then we present our main result in this section.

\begin{theorem}
\label{gmq}Under the \textbf{Assumptions} \textbf{1, 2} and \textbf{3}, RBSDE%
$(\xi ,f,L)$ admits a maximal bounded solution.
\end{theorem}

\proof%
First, notice that $(Y,Z,K)$ is the solution of RBSDE$(\xi ,f,L)$ if and
only if $(Y^{b},Z^{b},K^{b})$ is the solution of the RBSDE$(\xi
^{b},f^{b},L^{b})$, where
\[
(Y^{b},Z^{b},K^{b})=(Y-b,Z,K),
\]
and
\[
(\xi ^{b},f^{b}(t,y,z),L^{b})=(\xi -b,f(s,y+b,z),L-b).
\]
Notice that $(\xi ^{b},f^{b},L^{b})$ satisfies \textbf{Assumption 1, 2 }and $%
-2b\leq L^{b}\leq 0$. So in the following, we assume that the barrier $L$ is
a negative bounded process.

For $C>0$, set $g^{C}:\mathbb{R\rightarrow R}$ be a continuous
function, such that $0\leq g^{C}(y)\leq 1$, $\forall y\in
\mathbb{R}$, and
\begin{eqnarray}
g^{C}(y) &=&1\mbox{, if }\left| y\right| \leq C,  \label{gfun} \\
g^{C}(y) &=&0\mbox{, if }\left| y\right| \geq 2C.  \nonumber
\end{eqnarray}
Denote $f^{C}(t,y,z)=g^{C}(y)f(t,y,z)$; then
\begin{eqnarray*}
\left| f^{C}(t,y,z)\right| &\leq &g^{C}(y)(\varphi (\left| y\right|
)+A\left| z\right| ^{2}) \\
&\leq &1_{[-2C,2C]}(y)(\varphi (\left| y\right| )+A\left| z\right| ^{2}) \\
&\leq &\varphi (2C)+A\left| z\right| ^{2}.
\end{eqnarray*}
From the theorem 1 in \cite{KLQT}, there exists a maximal solution $%
(Y^{C},Z^{C},K^{C})$ to the RBSDE$(\xi ,f^{C},L)$%
\begin{eqnarray}
Y_{t}^{C} &=&\xi
+\int_{t}^{T}g^{C}(Y_{s}^{C})f(s,Y_{s}^{C},Z_{s}^{C})ds-%
\int_{t}^{T}Z_{s}^{C}dB_{s}+K_{T}^{C}-K_{t}^{C},  \label{RBSDE-qb1} \\
Y_{t}^{C} &\geq &L_{t},\int_{0}^{T}(Y_{t}^{C}-L_{t})dK_{t}^{C}=0,\mbox{ a.e..%
}  \nonumber
\end{eqnarray}
We choose $n\geq 2$ even, and $a\in \mathbb{R}$; applying
It\^{o}'s formula to $e^{at}(Y_{t}^{C})^{n}$, we have
\begin{eqnarray}
e^{at}(Y_{t}^{C})^{n} &=&e^{aT}\xi
^{n}+n%
\int_{t}^{T}e^{as}(Y_{s}^{C})^{n-1}g^{C}(Y_{s}^{C})f(s,Y_{s}^{C},Z_{s}^{C})ds-n\int_{t}^{T}e^{as}(Y_{s}^{C})^{n-1}Z_{s}^{C}dB_{s}
\label{est-qb1} \\
&&-\frac{n(n-1)}{2}\int_{t}^{T}e^{as}(Y_{s}^{C})^{n-2}\left|
Z_{s}^{C}\right|
^{2}ds+n\int_{t}^{T}e^{as}(Y_{s}^{C})^{n-1}dK_{s}^{C}-a%
\int_{t}^{T}e^{as}(Y_{s}^{C})^{n}ds.  \nonumber
\end{eqnarray}
From \textbf{Assumption 2} and the fact that $n$ is even, we have
\begin{eqnarray*}
yf(s,y,z) &\leq &yf(s,0,z)+\mu y^{2}, \\
y^{n-1}f(s,y,z) &\leq &y^{n-1}f(s,0,z)+\mu y^{n}.
\end{eqnarray*}
With $0\leq g^{C}(y)\leq 1$, we get
\begin{eqnarray*}
g^{C}(y)y^{n-1}f(s,y,z) &\leq &g^{C}(y)\left| y\right| ^{n-1}f(s,0,z)+\mu
y^{n} \\
&\leq &g^{C}(y)\left| y\right| ^{n-1}(\varphi (0)+A\left| z\right| ^{2})+\mu
y^{n} \\
&\leq &(\frac{1}{n}+\frac{n-1}{n}\left| y\right| ^{n})\varphi (0)+A\left|
z\right| ^{2}g^{C}(y)\left| y\right| ^{n-1}+\mu y^{n} \\
&\leq &(1+y^{n})\varphi (0)+2CA\left| z\right| ^{2}y^{n-2}+\mu y^{n}.
\end{eqnarray*}
Substitute it into (\ref{est-qb1}), then
\begin{eqnarray*}
e^{at}(Y_{t}^{C})^{n} &\leq &e^{aT}\xi ^{n}+\frac{n\varphi (0)}{a}%
(e^{aT}-e^{at})+(n\varphi (0)+n\mu -a)\int_{t}^{T}e^{as}(Y_{s}^{C})^{n}ds \\
&&+(2nCA-\frac{n(n-1)}{2})\int_{t}^{T}e^{as}(Y_{s}^{C})^{n-2}\left|
Z_{s}^{C}\right| ^{2}ds+n\int_{t}^{T}e^{as}(L_{s})^{n-1}dK_{s}^{C} \\
&&-n\int_{t}^{T}e^{as}(Y_{s}^{C})^{n-1}Z_{s}^{C}dB_{s}.
\end{eqnarray*}
Notice that since $K^{C}$ is an increasing process, $n$ is even and $L\leq 0$%
, we get immediately
\[
\int_{t}^{T}e^{as}(L_{s})^{n-1}dK_{s}^{C}\leq 0.
\]
If we choose $n$ and $a$ satisfying
\[
n-1\geq 4CA,a=n(\varphi (0)+\mu ),
\]
then
\[
e^{at}(Y_{t}^{C})^{n}\leq e^{aT}\xi ^{n}+\frac{n\varphi (0)}{a}%
(e^{aT}-e^{at})-n\int_{t}^{T}e^{as}(Y_{s}^{C})^{n-1}Z_{s}^{C}dB_{s}.
\]
It follows that
\[
e^{at}(Y_{t}^{C})^{n}\leq E[e^{aT}(\xi ^{n}+\frac{n\varphi (0)}{a})|\mathcal{%
F}_{t}]\leq e^{aT}(\left\| \xi \right\| _{\infty }^{n}+1),
\]
at last we get
\[
(Y_{t}^{C})^{n}\leq e^{a(T-t)}(\left\| \xi \right\| _{\infty }^{n}+1)\leq
(e^{aT}\vee 1)(\left\| \xi \right\| _{\infty }^{n}+1).
\]
Since $a=n(\varphi (0)+\mu )$, it follows that
\[
\left| Y_{t}^{C}\right| \leq (e^{(\varphi (0)+\mu )T}\vee 1)(\left\| \xi
\right\| _{\infty }^{n}+1)^{\frac{1}{n}}\leq (e^{(\varphi (0)+\mu )T}\vee
1)(\left\| \xi \right\| _{\infty }+1).
\]

If $C$ is chosen to satisfy $C\geq (e^{(\varphi (0)+\mu )T}\vee 1)(\left\|
\xi \right\| _{\infty }+1)$, then we have $\left| Y_{t}^{C}\right| \leq C$,
which implies $g^{C}(Y_{t}^{C})=1$, for $0\leq t\leq T$. So, $%
(Y^{C},Z^{C},K^{C})$ is the solution of the RBSDE$(\xi ,f,L)$. $\square $

\section{The case $f(t,y,z)=\left| z\right| ^{2}$}

In this section we consider the case $f(t,y,z)=\left| z\right| ^{2}$, which
corresponds to the RBSDE
\begin{eqnarray}
Y_{t} &=&\xi +\int_{t}^{T}\left| Z_{s}\right|
^{2}ds+K_{T}-K_{t}-\int_{t}^{T}Z_{s}dB_{s},  \label{RBSDE-p21} \\
Y_{t} &\geq &L_{t},\int_{0}^{T}(Y_{t}-L_{t})dK_{t}=0.  \nonumber
\end{eqnarray}
Then we have

\begin{theorem}
Under the assumption $E(\sup_{0\leq t\leq T}e^{2L_{t}})<+\infty $, the RBSDE$%
(\xi ,f,L)$ (\ref{RBSDE-p21}) admits a solution if and only if $E(e^{2\xi
})<+\infty $.
\end{theorem}

\proof%
For the necessary part, let $(Y,Z,K)$ be a solution of the RBSDE (\ref
{RBSDE-p21}). By It\^{o}'s formula, we get
\begin{eqnarray}
e^{2Y_{t}} &=&e^{2\xi
}+2\int_{t}^{T}e^{2Y_{s}}dK_{s}-2\int_{t}^{T}e^{Y_{s}}Z_{s}dB_{s}
\label{RBSDE-p22} \\
&=&e^{2Y_{0}}+2\int_{0}^{t}e^{2Y_{s}}Z_{s}dB_{s}-2%
\int_{0}^{t}e^{2Y_{s}}dK_{s}.  \nonumber
\end{eqnarray}
Let for all $n$, $\tau _{n}=\inf \{t:Y_{t}\geq n\}\wedge T$, then $%
M_{t\wedge \tau _{n}}=2\int_{0}^{t\wedge \tau _{n}}e^{2Y_{s}}Z_{s}dB_{s}$ is
a martingale, and we have
\[
E[e^{2Y_{\tau _{n}}}]=E[e^{2Y_{0}}-2\int_{0}^{t}e^{2Y_{s}}dK_{s}]\leq
E[e^{2Y_{0}}],
\]
in view of $2\int_{0}^{t}e^{2Y_{s}}dK_{s}\geq 0$. Finally, since $\tau
_{n}\nearrow T$, when $n\rightarrow \infty $:
\[
E[\underline{\lim }_{n\rightarrow \infty }e^{2Y_{\tau _{n}}}]=E[e^{2\xi
}]\leq E[e^{2Y_{0}}]<\infty ,
\]
follows from Fatou's Lemma.

Now we suppose $E(e^{2\xi })<+\infty $, set $\widetilde{L}%
_{t}=L_{t}1_{\{t<T\}}+\xi 1_{\{t=T\}}$ and
\[
N_{t}=S_{t}(e^{2\widetilde{L}})=ess\sup_{\tau \in \mathcal{T}_{t,T}}E[e^{2%
\widetilde{L}_{\tau }}|\mathcal{F}_{t}],
\]
where $S_{t}(\eta )$ denotes the Snell envelope of $\eta $ (See El Karoui
\cite{E79}), $\mathcal{T}_{t,T}$ is the set of all stopping times valued in $%
[t,T]$. Since
\[
E[\sup_{0\leq t\leq T}e^{2\widetilde{L}_{t}}]\leq E[\sup_{0\leq t\leq
T}e^{2L_{t}}+e^{2\xi }]<+\infty ,
\]
using the results of Snell envelope, we know that $N$ is a supermartingale,
so it admits the following decomposition: for an increasing integrable
process $\overline{K}$,
\[
N_{t}=N_{0}+\int_{0}^{t}\overline{Z}_{s}dB_{s}-\overline{K}_{t}.
\]
Applying It\^{o}'s formula to $\log N_{t}$, we get
\[
\frac{1}{2}\log N_{t}=\frac{1}{2}\log N_{0}+\frac{1}{2}\int_{0}^{t}\frac{%
\overline{Z}_{s}}{N_{s}}dB_{s}-\frac{1}{4}\int_{0}^{t}(\frac{\overline{Z}_{s}%
}{N_{s}})^{2}ds-\frac{1}{2}\int_{0}^{t}\frac{1}{N_{s}}d\overline{K}_{s}.
\]
Set $Y_{t}=\frac{1}{2}\log N_{t}$, $Z_{t}=\frac{\overline{Z}_{t}}{2N_{t}}$, $%
K_{t}=\frac{1}{2}\int_{0}^{t}\frac{1}{N_{s}}d\overline{K}_{s}$, then the
triple satisfies
\begin{equation}
Y_{t}=\xi +\int_{t}^{T}Z_{s}^{2}ds+K_{T}-K_{t}-\int_{t}^{T}Z_{s}dB_{s}.
\label{RBSDE-qb-s}
\end{equation}
Thanks to the results on the Snell envelope, we know that $N_{t}\geq e^{2%
\widetilde{L}_{t}}$ and $\int_{0}^{T}(N_{t}-e^{2\widetilde{L}_{t}})d%
\overline{K}_{t}=0$. The first implies
\[
Y_{t}\geq \widetilde{L}_{t}\geq L_{t}.
\]
Obviously, $N_{t}>0$, $0\leq t\leq T$, so $K$ is increasing. Consider the
stopping time $D_{t}:=\inf \{t\leq u\leq T;Y_{u}=L_{u}\}\wedge T$, then it
satisfies $D_{t}=\inf \{t\leq u\leq T;N_{u}=e^{2L_{u}}\}\wedge T$. By the
continuity of $\overline{K}$, we get $\overline{K}_{D_{t}}-\overline{K}_{t}=0
$, which implies $K_{D_{t}}-K_{t}=0$. It follows that
\[
\int_{0}^{T}(Y_{t}-L_{t})dK_{t}=0.
\]
Now the rest is to prove $Y_{t}\in \mathbf{S}^{2}(0,T)$, $Z_{t}\in \mathbf{H}%
_{d}^{2}(0,T)$, and $K_{t}\in \mathbf{A}^{2}(0,T)$. With Jensen's inequality
\begin{eqnarray*}
Y_{t} &=&\frac{1}{2}\log N_{t}=\frac{1}{2}\log [ess\sup_{\tau \in \mathcal{T}%
_{t,T}}E[e^{2\widetilde{L}_{\tau }}|\mathcal{F}_{t}]] \\
&\geq &\frac{1}{2}\log [\exp (ess\sup_{\tau \in \mathcal{T}_{t,T}}E[2%
\widetilde{L}_{\tau }|\mathcal{F}_{t}])] \\
&=&ess\sup_{\tau \in \mathcal{T}_{t,T}}E[\widetilde{L}_{\tau }|\mathcal{F}%
_{t}]\geq E[\xi |\mathcal{F}_{t}]\geq U_{t},
\end{eqnarray*}
where $U_{t}=-E[\xi ^{-}|\mathcal{F}_{t}]$. For all $a>0$, define
\[
\tau _{a}=\inf \{t;\left| N_{t}\right| >a,\int_{0}^{t}(\frac{\overline{Z}_{s}%
}{N_{s}})^{2}ds>a,\left| \int_{0}^{t}\frac{\overline{Z}_{s}}{N_{s}}%
dB_{s}\right| >a\}.
\]
From (\ref{RBSDE-qb-s}), we get for $0\leq t\leq T$%
\begin{eqnarray*}
0 &\leq &\int_{0}^{t}Z_{s}^{2}ds=Y_{0}-Y_{t}+\int_{0}^{t}Z_{s}dB_{s}-K_{t} \\
&\leq &Y_{0}-U_{t}+\int_{0}^{t}Z_{s}dB_{s}.
\end{eqnarray*}
Then
\[
(\int_{0}^{\tau _{a}}Z_{s}^{2}ds)^{2}\leq 3(Y_{0})^{2}+3(U_{\tau
_{a}})^{2}+3(\int_{0}^{\tau _{a}}Z_{s}dB_{s})^{2}.
\]
Taking the expectation, using the Jensen's inequality and $3x\leq \frac{x^{2}%
}{2}+\frac{9}{2}$, we obtain
\begin{eqnarray*}
E(\int_{0}^{\tau _{a}}Z_{s}^{2}ds)^{2} &\leq &\frac{3}{4}(\log
N_{0})^{2}+3E(\xi ^{-})^{2}+\frac{1}{2}(E(\int_{0}^{\tau
_{a}}Z_{s}^{2}ds))^{2}+\frac{9}{2} \\
&\leq &\frac{3}{4}(\log N_{0})^{2}+3E(\xi ^{-})^{2}+\frac{1}{2}%
E(\int_{0}^{\tau _{a}}Z_{s}^{2}ds)^{2}+\frac{9}{2},
\end{eqnarray*}
so
\[
E(\int_{0}^{\tau _{a}}Z_{s}^{2}ds)^{2}\leq \frac{3}{2}(\log
N_{0})^{2}+6E(\xi ^{-})^{2}+9\leq C.
\]
Since $\tau _{a}\nearrow T$ when $a\rightarrow +\infty $, we get to the
limit, and with the Schwartz inequality
\[
E\int_{0}^{T}Z_{s}^{2}ds\leq (E(\int_{0}^{T}Z_{s}^{2}ds)^{2})^{\frac{1}{2}%
}\leq C.
\]
So $Z\in \mathbf{H}_{d}^{2}(0,T)$. From (\ref{RBSDE-qb-s}), we get for $%
0\leq t\leq T$%
\begin{eqnarray*}
0 &\leq &K_{t}=Y_{0}-Y_{t}+\int_{0}^{t}Z_{s}dB_{s}-\int_{0}^{t}Z_{s}^{2}ds \\
&\leq &Y_{0}-Y_{t}+\int_{0}^{t}Z_{s}dB_{s}.
\end{eqnarray*}
Notice that $K$ is increasing, so it's sufficient to prove $%
E[K_{T}^{2}]<+\infty $. Squaring the inequality on both sides and taking
expectation, we obtain
\[
E[(K_{T})^{2}]\leq 3Y_{0}^{2}+3E[\xi ^{2}]+3E\int_{0}^{T}Z_{s}^{2}ds\leq C.
\]
We consider now $Y$; again from (\ref{RBSDE-qb-s}),
\[
Y_{t}=Y_{0}-K_{t}+\int_{0}^{t}Z_{s}dB_{s}-\int_{0}^{t}Z_{s}^{2}ds,
\]
so
\[
(Y_{t})^{2}\leq 4\left( Y_{0}\right) ^{2}+4\left( K_{t}\right) ^{2}+4\left(
\int_{0}^{t}Z_{s}dB_{s}\right) ^{2}+4\left( \int_{0}^{t}Z_{s}^{2}ds\right)
^{2}.
\]
Then by the Bukholder-Davis-Gundy inequality, we get
\begin{eqnarray*}
E[\sup_{0\leq t\leq T}(Y_{t})^{2}] &\leq &4\left( Y_{0}\right)
^{2}+4E[K_{T}^{2}]+4E[\sup_{0\leq t\leq T}\left(
\int_{0}^{t}Z_{s}dB_{s}\right) ^{2}]+4E\left( \int_{0}^{T}Z_{s}^{2}ds\right)
^{2} \\
&\leq &4\left( Y_{0}\right) ^{2}+4E[K_{T}^{2}]+CE\left(
\int_{0}^{t}Z_{s}^{2}dB_{s}\right) +4E\left( \int_{0}^{T}Z_{s}^{2}ds\right)
^{2}\leq C,
\end{eqnarray*}
i.e. $Y\in \mathbf{S}^{2}(0,T)$. $\square $

\section{The case when $f$ is linear increasing in $z$}

In this section, we assume that the coefficient $f$ satisfies

\textbf{Assumption 6. }(i) $f(\cdot ,y,z)$ is progressively measurable, and $%
E\int_{0}^{T}f^{2}(t,0,0)dt<+\infty $;

(ii) for $\mu \in \mathbb{R}$, $\forall (t,z)\in [0,T]\times
\mathbb{R}^{d}$
and $y,y^{\prime }\in \mathbb{R,}$%
\[
(y-y^{\prime })(f(t,y,z)-f(t,y^{\prime },z))\leq \mu (y-y^{\prime })^{2};
\]

(iii) there exists a nonegative, continuous, increasing function $\varphi :%
\mathbb{R}^{+}\rightarrow \mathbb{R}^{+}$, with $\varphi (0)=0$, s.t. $%
\forall (t,y,z)\in [0,T]\times \mathbb{R\times R}^{d},$%
\[
\left| f(t,y,z)\right| \leq \left| g_{t}\right| +\varphi (\left| y\right|
)+\beta \left| z\right| ,
\]
where $g_{t}\in \mathbf{H}^{2}(0,T)$;

(iv) for $t\in [0,T]$, $(y,z)\rightarrow f(t,y,z)$ is continuous.

If $\varphi (x)=\left| x\right| $, then $f$ is linear increasing in $y$ and $%
z$. Matoussi proved in \cite{M97} that when $\xi \in \mathbf{L}^{2}(\mathcal{%
F}_{T})$ and $L\in \mathbf{S}^{2}(0,T)$, there exists a triple $(Y,Z,K)$
which is solution of the RBSDE$(\xi ,f,L)$.

Our result of this section is the following:

\begin{theorem}
\label{exist-zl}Suppose that $\xi \in \mathbf{L}^{2}(\mathcal{F}_{T})$, $f$
and $L$ satisfy \textbf{Assumption} \textbf{6} and \textbf{3}, respectively,
then the RBSDE$(\xi ,f,L)$ has a minimal solution $(Y,Z,K)\in \mathbf{S}%
^{2}(0,T)\times \mathbf{H}_{d}^{2}(0,T)\times \mathbf{A}^{2}(0,T)$, which
satisfies
\[
Y_{t}=\xi
+\int_{t}^{T}f(s,Y_{s},Z_{s})ds+K_{T}-K_{t}-\int_{t}^{T}Z_{s}dB_{s},
\]
$Y_{t}\geq L_{t}$, and $\int_{0}^{T}(Y_{s}-L_{s})dK_{s}=0$.
\end{theorem}

First we note that the triple $(Y,Z,K)$ solves the RBSDE$(\xi ,f,L)$, if and
only if the triple
\begin{equation}
(\overline{Y}_{t},\overline{Z}_{t},\overline{K}_{t}):=(e^{\lambda
t}Y_{t},e^{\lambda t}Z_{t},\int_{0}^{t}e^{\lambda s}dK_{s})  \label{trans}
\end{equation}
solves the RBSDE$(\overline{\xi },\overline{f},\overline{L})$, where
\[
(\overline{\xi },\overline{f}(t,y,z),\overline{L}_{t})=(\xi e^{\lambda
T},e^{\lambda t}f(t,e^{-\lambda t}y,e^{-\lambda t}z)-\lambda y,e^{\lambda
t}L_{t}).
\]

If we choose $\lambda =\mu $, then the coefficient $\overline{f}$ satisfies
the same assumptions as in \textbf{Assumption 6}, with (ii) replaced by

(ii') $(y-y^{\prime })(f(t,y,z)-f(t,y^{\prime },z))\leq 0$.

Since we are in the $1$-dimensional case, (ii') means that $f$ is decreasing
on $y$. From another part $\overline{\xi }$ still belongs to $\mathbf{L}^{2}(%
\mathcal{F}_{T})$ and the barrier $\overline{L}$ still satisfy the
assumptions \textbf{Assumption} \textbf{3}. So\textbf{\ }in the following,
we shall work under \textbf{Assumption 6' }with (ii) replaced by (ii').

Before proving this theorem, we consider an estimate result and a monotonic
stability theorem for RBSDEs.

\begin{lemma}
\label{est-res}We consider RBSDE$(\xi ,g,L)$, with $\xi \in \mathbf{L}^{2}(%
\mathcal{F}_{T})$, $g$ and $L$ satisfy \textbf{Assumption} \textbf{6'} and
\textbf{3}. Moreover $g(t,y,z)$ is Lipschitz in $z$. Then we have the
following estimation
\begin{eqnarray*}
&&E[\sup_{0\leq t\leq T}\left| y_{t}\right| ^{2}+\int_{0}^{T}\left|
z_{s}\right| ds+\left| k_{T}\right| ^{2}] \\
&\leq &C_{\beta }E[\left| \xi \right| ^{2}+\int_{0}^{T}g_{s}^{2}ds+\varphi
^{2}(b)+\varphi ^{2}(2T)+1]
\end{eqnarray*}
where $(y_{t},z_{t},k_{t})_{0\leq t\leq T}$ is the solution of RBSDE$(\xi
,g,L)$. $C_{\beta }$ is a constant only depends on $\beta $, $T$ and $b$.
\end{lemma}

\begin{remark}
The constant $C_{\beta }$ does not depend on Lipschitz coefficient of $g$ on
$z$.
\end{remark}

\proof%
Since $g\,$is Lipschitz in $z$, by the theorem 2 in \cite{LMX}, the RBSDE$%
(\xi ,g,L)$ admits the unique solution $(y_{t},z_{t},k_{t})_{0\leq t\leq T}$%
. Apply It\^{o}'s formula to $\left| y_{t}\right| ^{2}$, in view of $%
yg(t,y,z)\leq g(t,0,0)\left| y\right| +\beta \left| y\right| \left| z\right|
$ and $\sup_{0\leq t\leq T}\left| L_{t}\right| \leq b$, we get
\begin{eqnarray*}
E[\left| y_{t}\right| ^{2}+\int_{t}^{T}\left| z_{s}\right| ^{2}ds]
&=&E[\left| \xi \right|
^{2}+2\int_{t}^{T}y_{s}g(s,y_{s},z_{s})ds+2\int_{t}^{T}L_{s}dk_{s}] \\
&\leq &E[\left| \xi \right| ^{2}+2\int_{t}^{T}y_{s}g_{s}ds+2\beta
\int_{t}^{T}y_{s}z_{s}ds+2b(k_{T}-k_{t})].
\end{eqnarray*}
It follows that
\begin{eqnarray*}
&&E[\left| y_{t}\right| ^{2}+\frac{1}{2}\int_{t}^{T}\left| z_{s}\right|
^{2}ds] \\
&\leq &E[\left| \xi \right| ^{2}+\int_{t}^{T}g_{s}^{2}ds+(1+2\beta
^{2})\int_{t}^{T}\left| y_{s}\right| ^{2}ds+2b(k_{T}-k_{t})].
\end{eqnarray*}
By Gronwall's inequality, we know there exists a constant $c_{1}$ depending
on $\beta $ and $T$, such that for $t\in [0,T]$,
\begin{equation}
E[\left| y_{t}\right| ^{2}]\leq c_{1}E[\left| \xi \right|
^{2}+\int_{0}^{T}g_{s}^{2}ds+b(k_{T}-k_{t})].  \label{est-y1}
\end{equation}
It follows that
\begin{equation}
E[\int_{t}^{T}\left| z_{s}\right| ^{2}ds]\leq 2(1+(1+2\beta
^{2})T)c_{1}E[\left| \xi \right|
^{2}+\int_{0}^{T}g_{s}^{2}ds+b(k_{T}-k_{t})].  \label{est-z1}
\end{equation}

Now we estimate the increasing process $k$ by approximation. Take $z$ as a
known process, without losing of generality, we write $g(t,y)$ for $%
g(t,y,z_{t})$, here $g(t,0)=g(t,0,z_{t})$ is a process in $\mathbf{H}%
^{2}(0,T)$ in view of linear increasing property of $g$ on $z$.

For $m$, $p\in \mathbb{N}$, set $\xi ^{m,p}=(\xi \vee (-p))\wedge m$, $%
g^{m,p}(t,u)=g(t,u)-g_{t}+(g_{t}\vee (-p))\wedge m$. We consider RBSDE$(\xi
^{m,p},g^{m,p},L)$,
\begin{eqnarray}
y_{t}^{m,p} &=&\xi
^{m,p}+\int_{t}^{T}g^{m,p}(s,y_{s}^{m,p})ds+k_{T}^{m,p}-k_{t}^{m,p}-%
\int_{t}^{T}z_{s}^{m,p}dB_{s},  \label{appro} \\
y_{t}^{m,p} &\geq &L_{t},\int_{0}^{T}(y_{t}^{m,p}-L_{t})dk_{t}^{m,p}=0.
\nonumber
\end{eqnarray}
It is easy to check that $(y^{m,p},z^{m,p},k^{m,p})$ is the solution of RBSDE%
$(\xi ^{m,p},g^{m,p},L)$, if and only if $(\widehat{y}^{m,p},\widehat{z}%
^{m,p},\widehat{k}^{m,p})$ is the solution of RBSDE$(\widehat{\xi }^{m,p},%
\widehat{g}^{m,p},\widehat{L})$, where
\[
(\widehat{y}_{t}^{m,p},\widehat{z}_{t}^{m,p},\widehat{k}%
_{t}^{m,p})=(y_{t}^{m,p}+m(t-2(T\vee 1)),z_{t}^{m,p},k_{t}^{m,p}),
\]
and
\begin{eqnarray*}
\widehat{\xi }^{m,p} &=&\xi ^{m,p}+mT-2m(T\vee 1), \\
\widehat{g}^{m,p}(t,y) &=&g^{m,p}(t,y-m(t-2(T\vee 1)))-m, \\
\widehat{L}_{t} &=&L_{t}+m(t-2(T\vee 1)).
\end{eqnarray*}
Without losing of generality, we set $T\geq 1$. Since $\xi ^{m,p}$ and $%
g_{t}^{m,p}\leq m$, we have $\widehat{\xi }^{m,p}$ and $\widehat{g}%
_{t}^{m,p}\leq 0$. By (\ref{appro}),
\[
\widehat{k}_{T}^{m,p}-\widehat{k}_{t}^{m,p}=\widehat{y}_{t}^{m,p}-\widehat{%
\xi }^{m,p}-\int_{t}^{T}\widehat{g}^{m,p}(s,\widehat{y}_{s}^{m,p})ds+%
\int_{t}^{T}\widehat{z}_{s}^{m,p}dB_{s},
\]
taking square and expectation on the both sides, we get
\begin{equation}
E[(\widehat{k}_{T}^{m,p}-\widehat{k}_{t}^{m,p})^{2}]\leq 4E[(\widehat{y}%
_{t}^{m,p})^{2}+(\widehat{\xi }^{m,p})^{2}+(\int_{t}^{T}\widehat{g}^{m,p}(s,%
\widehat{y}_{s}^{m,p})ds)^{2}+\int_{t}^{T}\left| \widehat{z}%
_{s}^{m,p}\right| ^{2}ds].  \label{est-k1}
\end{equation}

In order to estimate the first and the last form on the left side, we apply
It\^{o}'s formula to $\left| \widehat{y}_{t}^{m,p}\right| ^{2}$, and get the
following with Gronwall inequality,
\begin{eqnarray}
&&E[\left| \widehat{y}_{t}^{m,p}\right| ^{2}+\int_{t}^{T}\left| \widehat{z}%
_{s}^{m,p}\right| ^{2}ds]  \label{est-yz1} \\
&\leq &c_{2}E[\left| \widehat{\xi }^{m,p}\right| ^{2}+\int_{t}^{T}(\widehat{g%
}_{s}^{m,p})^{2}ds+\int_{t}^{T}L_{s}d\widehat{k}_{s}^{m,p}],  \nonumber
\end{eqnarray}
where $c_{2}$ is a constant only depends on $T$. For the third term, let us
recall a comparison result of $\widehat{y}_{t}^{m,p}$ in step 2 of the proof
of theorem 2 in \cite{LMX},
\[
\widetilde{y}_{t}^{m,p}\leq \widehat{y}_{t}^{m,p}\leq \overline{y}%
_{t}^{m,p},
\]
where $\widetilde{y}_{t}^{m,p}$ is the solution of BSDE$(\widehat{\xi }%
^{m,p},\widehat{g}^{m,p})$, i.e.
\begin{equation}
\widetilde{y}_{t}^{m,p}=\widehat{\xi }^{m,p}+\int_{t}^{T}\widehat{g}^{m,p}(s,%
\widetilde{y}_{s}^{m,p})ds-\int_{t}^{T}\widetilde{z}_{s}^{m,p}dB_{s},
\label{est-bsde}
\end{equation}
and
\[
\overline{y}_{t}^{m,p}=ess\sup_{\tau \in \mathcal{T}_{t,T}}E[(\widehat{L}%
_{\tau })^{+}1_{\{\tau <T\}}+(\widehat{\xi }^{m,p})^{+}1_{\{\tau =T\}}|%
\mathcal{F}_{t}],
\]
where $\mathcal{T}_{t,T}$ is the set of stoppng times valued in $[t,T]$.
Moreover, we have $\sup_{0\leq s\leq T}\overline{y}_{s}^{m,p}=\sup_{0\leq
s\leq T}\widehat{L}_{s}$.

Since $\widehat{g}^{m,p}$ is decreasing in $y$, we get
\[
\widehat{g}^{m,p}(s,\overline{y}_{t}^{m,p})\leq \widehat{g}^{m,p}(s,\widehat{%
y}_{s}^{m,p})\leq \widehat{g}^{m,p}(s,\widetilde{y}_{s}^{m,p}).
\]
So to estimate $E[(\int_{t}^{T}\widehat{g}^{m,p}(s,\widehat{y}%
_{s}^{m,p})ds)^{2}]$, it is sufficante to get the estimations of $%
E[(\int_{t}^{T}\widehat{g}^{m,p}(s,\widetilde{y}_{s}^{m,p})ds)^{2}]$ and $%
E[(\int_{t}^{T}\widehat{g}^{m,p}(s,\overline{y}_{s}^{m,p})ds)^{2}]$. First
we know that
\begin{eqnarray}
E[(\int_{t}^{T}\widehat{g}^{m,p}(s,\widetilde{y}_{s}^{m,p})ds)^{2}] &\leq
&3E[\left| \widehat{\xi }^{m,p}\right| ^{2}+\left| \widetilde{y}%
_{t}^{m,p}\right| ^{2}+\int_{t}^{T}\left| \widetilde{z}_{s}^{m,p}\right|
^{2}ds]  \label{est-g1} \\
&\leq &c_{3}E[\left| \widehat{\xi }^{m,p}\right| ^{2}+\int_{t}^{T}(\widehat{g%
}_{s}^{m,p})^{2}ds],  \nonumber
\end{eqnarray}
in view of estimate result of BSDE(\ref{est-bsde}). Here $c_{3}$ is a
constant only depends on $T$. Then with the presentation of $\overline{y}%
_{t}^{m,p}$, we have
\begin{equation}
E[(\int_{t}^{T}\widehat{g}^{m,p}(s,\overline{y}_{s}^{m,p})ds)^{2}]\leq
E[2T\int_{0}^{T}(\widehat{g}_{s}^{m,p})^{2}ds+2T\varphi ^{2}(\sup_{0\leq
t\leq T}(\widehat{L}_{t})^{+})]  \label{est-g2}
\end{equation}
From (\ref{est-k1}), with (\ref{est-yz1}), (\ref{est-g1}) and (\ref{est-g2}%
), we have
\begin{eqnarray*}
E[(k_{T}^{m,p}-k_{t}^{m,p})^{2}] &=&E[(\widehat{k}_{T}^{m,p}-\widehat{k}%
_{t}^{m,p})^{2}] \\
&\leq &c_{4}E[\left| \widehat{\xi }^{m,p}\right| ^{2}+\int_{t}^{T}(\widehat{g%
}_{s}^{m,p})^{2}ds+\int_{t}^{T}L_{s}d\widehat{k}_{s}^{m,p}+\varphi
^{2}(\sup_{0\leq t\leq T}(\widehat{L}_{t})^{+})] \\
&\leq &c_{4}E[2\left| \xi ^{m,p}\right|
^{2}+4\int_{t}^{T}(g_{s}^{m,p})^{2}ds+2c_{4}b^{2}+\varphi ^{2}(b)] \\
&&+4c_{4}(m^{2}T^{2}+\varphi ^{2}(2mT))+\frac{1}{2}%
E[(k_{T}^{m,p}-k_{t}^{m,p})^{2}],
\end{eqnarray*}
where $c_{4}=c_{2}\vee c_{3}\vee (2T)$, which only depends on $T$. It
follows that
\begin{eqnarray*}
&&E[(k_{T}^{m,p}-k_{t}^{m,p})^{2}] \\
&\leq &c_{5}E[\left| \xi ^{m,p}\right|
^{2}+\int_{t}^{T}(g_{s}^{m,p})^{2}ds+b^{2}+\varphi
^{2}(b)]+c_{5}(m^{2}T^{2}+\varphi ^{2}(2mT)) \\
&\leq &c_{5}E[\left| \xi \right| ^{2}+\int_{t}^{T}g_{s}^{2}ds+b^{2}+\varphi
^{2}(b)]+c_{5}(m^{2}T^{2}+\varphi ^{2}(2mT)),
\end{eqnarray*}
where $c_{5}=4c_{4}^{2}\vee 8c_{4}$.

Now we consider the RBSDE$(\xi ^{p},g^{p},L)$, where $\xi ^{p}=\xi \vee (-p)$%
, $g^{p}(t,u)=g(t,u)-g_{t}+g_{t}\vee (-p)$. Thanks to the convergence result
in \cite{LMX}, we know that
\[
(y^{m,p},z^{m,p},k^{m,p})\rightarrow (y^{p},z^{p},k^{p})\mbox{ in }\mathbf{S}%
^{2}(0,T)\times \mathbf{H}_{d}^{2}(0,T)\times \mathbf{A}^{2}(0,T),
\]
where $(y^{p},z^{p},k^{p})$ is the solution of RBSDE$(\xi ^{p},g^{p},L)$.
Moreover, we have $dk_{t}^{p}\leq dk_{t}^{1,p}$, by comparison theorem. So
\begin{eqnarray*}
E[(k_{T}^{p}-k_{t}^{p})^{2}] &\leq &E[(k_{T}^{1,p}-k_{t}^{1,p})^{2}] \\
&\leq &c_{5}E[\left| \xi \right| ^{2}+\int_{t}^{T}g_{s}^{2}ds+b^{2}+\varphi
^{2}(b)]+c_{5}(T^{2}+\varphi ^{2}(2T)).
\end{eqnarray*}
Then let $p\rightarrow \infty $, thanks to the convergence result in \cite
{LMX}, we know
\[
(y^{p},z^{p},k^{p})\rightarrow (y,z,k)\mbox{ in
}\mathbf{S}^{2}(0,T)\times \mathbf{H}_{d}^{2}(0,T)\times
\mathbf{A}^{2}(0,T).
\]
In view \textbf{Assumption 6}-(iii), it follows that
\begin{eqnarray*}
E[(k_{T}-k_{t})^{2}] &\leq &c_{5}E[\left| \xi \right|
^{2}+\int_{t}^{T}(g(s,0,z_{s}))^{2}ds+b^{2}+\varphi ^{2}(b)] \\
&\leq &c_{5}E[\left| \xi \right| ^{2}+2\int_{t}^{T}g_{s}^{2}ds+2\beta
^{2}\int_{t}^{T}\left| z_{s}\right| ^{2}ds+b^{2}+\varphi ^{2}(b)] \\
&&+c_{5}(T^{2}+\varphi ^{2}(2T)).
\end{eqnarray*}
With (\ref{est-z1}), setting $c_{6}=c_{5}\vee (4\beta ^{2}(1+(1+2\beta
^{2})T)c_{1}+2)\vee c_{5}(b^{2}+T^{2})$, we get
\begin{eqnarray*}
E[(k_{T}-k_{t})^{2}] &\leq &c_{6}E[\left| \xi \right|
^{2}+\int_{t}^{T}g_{s}^{2}ds+b(k_{T}-k_{t})+\varphi ^{2}(b)+\varphi
^{2}(2T)+1] \\
&\leq &c_{6}E[\left| \xi \right|
^{2}+\int_{t}^{T}g_{s}^{2}ds+2c_{6}b^{2}+\varphi ^{2}(b)+\varphi ^{2}(2T)+1]
\\
&&+\frac{1}{2}E[(k_{T}-k_{t})^{2}]
\end{eqnarray*}
It follows that
\[
E[(k_{T}-k_{t})^{2}]\leq 2c_{6}E[\left| \xi \right|
^{2}+2\int_{t}^{T}g_{s}^{2}ds+\varphi ^{2}(b)+\varphi
^{2}(2T)+2c_{6}b^{2}+1].
\]
Consequantly, by (\ref{est-y1}) and (\ref{est-z1}), we obtain
\begin{eqnarray*}
&&\sup_{0\leq t\leq T}E[\left| y_{t}\right| ^{2}]+E[\int_{0}^{T}\left|
z_{s}\right| ds+\left| k_{T}\right| ^{2}] \\
&\leq &C_{\beta }E[\left| \xi \right| ^{2}+\int_{0}^{T}g_{s}^{2}ds+\varphi
^{2}(2T)+\varphi ^{2}(b)+1],
\end{eqnarray*}
where $C_{\beta }$ is a constant only depends on $\beta $, $T$ and $b$. The
final result follows from BDG inequality. $\square $

The proof of this theorem is step 1 and step 2 of the proof of theorem 4 in
\cite{KLQT}, with comparison theorem. So we omit it.

With these preparations, we begin our main proof.

\textbf{Proof of theorem \ref{exist-zl}. }The proof consists 4 step.

\textbf{Step 1. }Approximation. For $n\geq \beta $, we introduce the
following functions
\[
f_{n}(t,y,z)=\inf_{q\in \mathbf{Q}^{d}}\{f(t,y,q)+n\left| z-q\right| \},
\]
then we have

1. for all $(t,z)$, $y\rightarrow f_{n}(t,y,z)$ is non-increasing;

2. for all $(t,y)$, $z\rightarrow f_{n}(t,y,z)$ is $n$-Lipschitz;

3. for all $(t,y,z)$, $\left| f_{n}(t,y,z)\right| \leq \left| g_{t}\right|
+\varphi (\left| y\right| )+\beta \left| z\right| .$

Thanks to the results of \cite{LMX}, we know that for each $n\geq \beta $,
there exits a unique triple $(Y^{n},Z^{n},K^{n})$ satisfies the followings
\begin{eqnarray*}
Y_{t}^{n} &=&\xi
+\int_{t}^{T}f_{n}(s,Y_{s}^{n},Z_{s}^{n})ds+K_{T}^{n}-K_{t}^{n}-%
\int_{t}^{T}Z_{s}^{n}dB_{s}, \\
Y_{t}^{n} &\geq &L_{t},\int_{0}^{T}(Y_{t}^{n}-L_{t})dK_{t}^{n}=0.
\end{eqnarray*}

\textbf{Step 2}. Estimates results. Let $\alpha \geq 0$, be a real number to
be chosen later. We set
\[
U_{t}^{n}=e^{\alpha t}Y_{t}^{n},V_{t}^{n}=e^{\alpha
t}Z_{t}^{n},dJ_{t}^{n}=e^{\alpha t}dK_{t}^{n}.
\]
Then we know that $(U^{n},V^{n},J^{n})$ is the solution of the RBSDE
associated with $(\zeta ,F_{n},L^{\alpha })$, where
\[
\zeta =e^{\alpha T}\xi ,F_{n}(t,u,v)=e^{\alpha t}f_{n}(t,e^{-\alpha
t}u,e^{-\alpha t}v)-\alpha u,L_{t}^{\alpha }=e^{\alpha t}L_{t}.
\]
It is easy to check
\[
\left| F_{n}(t,u,v)\right| \leq e^{\alpha t}\left| g_{t}\right| +e^{\alpha
t}\varphi (\left| u\right| )+\alpha \left| u\right| +\beta \left| v\right| ,
\]
setting $\psi (u)=e^{\alpha T}\varphi (\left| u\right| )+\alpha \left|
u\right| $, with $\psi (u)=0$, we get that $F_{n}$ verifies\textbf{\
Assumption 6'}-(iii). Moreover
\begin{eqnarray*}
uF_{n}(t,u,v) &=&e^{\alpha t}uf_{n}(t,e^{-\alpha t}u,e^{-\alpha t}v)-\alpha
u^{2} \\
&\leq &ue^{\alpha t}g_{t}+\beta \left| u\right| \left| v\right| -\alpha
u^{2}.
\end{eqnarray*}
And $\sup_{0\leq t\leq T}L_{t}^{\alpha }\leq e^{\alpha T}\sup_{0\leq t\leq
T}L_{t}\leq e^{\alpha T}b$. Now we apply It\^{o} formula to $\left|
U^{n}\right| ^{2}$ on $[0,T]$, and get
\begin{eqnarray*}
&&\left| U_{t}^{n}\right| ^{2}+\int_{t}^{T}\left| V_{s}^{n}\right| ^{2}ds \\
&=&\left| \zeta \right|
^{2}+2\int_{t}^{T}U_{s}^{n}F_{n}(s,U_{s}^{n},V_{s}^{n})ds+2%
\int_{t}^{T}U_{s}^{n}dJ_{s}^{n}-2\int_{t}^{T}U_{s}^{n}V_{s}^{n}dB_{s} \\
&\leq &\left| \zeta \right| ^{2}+\int_{t}^{T}e^{2\alpha
s}g_{s}^{2}ds+(1+2\beta ^{2}-\alpha )\int_{t}^{T}\left| U_{s}^{n}\right|
^{2}ds+\frac{1}{2}\int_{t}^{T}\left| V_{s}^{n}\right| ds \\
&&+\theta e^{2\alpha T}b^{2}+\frac{1}{\theta }(J_{T}^{n}-J_{t}^{n})^{2}-2%
\int_{t}^{T}U_{s}^{n}V_{s}^{n}dB_{s},
\end{eqnarray*}
where $\theta $ is a constant to be decided later. By taking conditional
expectation, we get
\begin{eqnarray}
\left| U_{t}^{n}\right| ^{2}+\frac{1}{2}E[\int_{t}^{T}\left|
V_{s}^{n}\right| ^{2}ds|\mathcal{F}_{t}] &\leq &E[\left| \zeta \right|
^{2}+\int_{t}^{T}e^{2\alpha s}g_{s}^{2}ds+\theta e^{2\alpha T}b^{2}|\mathcal{%
F}_{t}]  \label{est-u1} \\
&&+(1+2\beta ^{2}-\alpha )E[\int_{t}^{T}\left| U_{s}^{n}\right| ^{2}ds|%
\mathcal{F}_{t}]+\frac{1}{\theta }E[(J_{T}^{n}-J_{t}^{n})^{2}|\mathcal{F}%
_{t}].  \nonumber
\end{eqnarray}
Since
\[
J_{T}^{n}-J_{t}^{n}=U_{t}^{n}-\zeta
-\int_{t}^{T}F_{n}(s,U_{s}^{n},V_{s}^{n})ds-\int_{t}^{T}V_{s}^{n}dB_{s},
\]
we have
\[
E[(J_{T}^{n}-J_{t}^{n})^{2}|\mathcal{F}_{t}]\leq 4\left| U_{t}^{n}\right|
^{2}+4E[\left| \zeta \right|
^{2}+(\int_{t}^{T}F_{n}(s,U_{s}^{n},V_{s}^{n})ds)^{2}+\int_{t}^{T}\left|
V_{s}^{n}\right| ^{2}ds|\mathcal{F}_{t}].
\]
Using the same approximation as in Lemma \ref{est-res}, except considering
conditional expectation $E[\cdot |\mathcal{F}_{t}]$ instead of expectation,
we deduce
\[
E[(J_{T}^{n}-J_{t}^{n})^{2}|\mathcal{F}_{t}]\leq c_{\beta }E[\left| \zeta
\right| ^{2}+\int_{t}^{T}e^{2\alpha s}g_{s}^{2}ds+\psi ^{2}(e^{\alpha
T}b)+\psi ^{2}(2T)+1|\mathcal{F}_{t}],
\]
where $c_{\beta }$ is a constant which only depends on $\beta $, $T$, $b$
and $\alpha $. Substitute it into (\ref{est-u1}), set $\alpha =1+2\beta ^{2}$%
, $\theta =c_{\beta }$, then we get,
\begin{eqnarray*}
\left| U_{t}^{n}\right| ^{2} &\leq &2E[\left| \zeta \right|
^{2}+\int_{t}^{T}F_{n}^{2}(s,0,0)ds|\mathcal{F}_{t}]+e^{\alpha T}(\varphi
(e^{\alpha T}b)+\varphi (2T)) \\
&&+\alpha (e^{\alpha T}b+2T)+1+c_{\beta }e^{2\alpha T}b^{2}.
\end{eqnarray*}
Recall the definition of $U^{n}$, we get
\begin{eqnarray*}
\left| Y_{t}^{n}\right| ^{2} &\leq &e^{-2\alpha t}(2E[e^{2\alpha T}\left|
\xi \right| ^{2}+\int_{t}^{T}e^{2\alpha s}g_{s}^{2}ds|\mathcal{F}%
_{t}]+e^{\alpha T}(\varphi (e^{\alpha T}b)+\varphi (2T)) \\
&&+\alpha (e^{\alpha T}b+2T)+1+c_{\beta }e^{2\alpha T}b^{2}).
\end{eqnarray*}
If we set $M_{t}=(e^{2\alpha T}2E[\left| \xi \right|
^{2}+\int_{t}^{T}g_{s}^{2}ds|\mathcal{F}_{t}]+e^{\alpha T}(\varphi
(e^{\alpha T}b)+\varphi (2T))+c_{\beta }e^{2\alpha T}b^{2}+\alpha (e^{\alpha
T}b+2T)+1)^{\frac{1}{2}}$, then
\begin{equation}
\left| Y_{t}^{n}\right| \leq M_{t},\forall t\in [0,T].  \label{est-y2}
\end{equation}

\textbf{Step 3}. Localisation.

First, we know that the sequence $(f_{n})_{n\geq \beta }$ is non-decreasing
in $n$, then from comparison theorem in \cite{LMX}, we get
\[
Y_{t}^{n}\leq Y_{t}^{n+1},\forall t\in [0,T],\forall n\geq \beta .
\]
Define $Y_{t}=\sup_{n\geq \beta }Y_{t}^{n}$.

We now consider the localisation procedure. For $m\in \mathbb{N}$, $m\geq b$%
, let $\tau _{m}$ be the following stopping time
\[
\tau _{m}=\inf \{t\in [0,T]:M_{t}+g_{t}\geq m\}\wedge T,
\]
and we introduce the stopped process $Y_{t}^{n,m}=Y_{t\wedge \tau _{m}}^{n}$%
, together with $Z_{t}^{n,m}=Z_{t}^{n}1_{\{t\leq \tau _{m}\}}$ and $%
K_{t}^{n,m}=K_{t\wedge \tau _{m}}^{n}$. Then $%
(Y_{t}^{n,m},Z_{t}^{n,m},K_{t}^{n,m})_{0\leq t\leq T}$ solved the following
RBSDE
\begin{eqnarray*}
Y_{t}^{n,m} &=&\xi ^{n,m}+\int_{t}^{T}1_{\{s\leq \tau
_{m}\}}f_{n}(s,Y_{s}^{n,m},Z_{s}^{n,m})ds+K_{T}^{n,m}-K_{t}^{n,m}-%
\int_{t}^{T}Z_{s}^{n,m}dB_{s}, \\
Y_{t}^{n,m} &\geq &L_{t},\int_{0}^{T}(Y_{t}^{n,m}-L_{t})dK_{t}^{n,m}=0.
\end{eqnarray*}
where $\xi ^{n,m}=Y_{\tau _{m}}^{n,m}=Y_{\tau _{m}}^{n}$.

Since $(Y^{n,m})_{n\geq \beta }$ is non-decreasing in $n$, with (\ref{est-y2}%
), we get $\sup_{n\geq \beta }\sup_{t\in [0,T]}\left| Y_{t}^{n,m}\right|
\leq m$. Set $\rho _{m}(y)=\frac{ym}{\max \{\left| y\right| ,m\}}$, then it
is easy to check that $(Y^{n,m},Z^{n,m},K^{n,m})$ verifies
\begin{eqnarray*}
Y_{t}^{n,m} &=&\xi ^{n,m}+\int_{t}^{T}1_{\{s\leq \tau _{m}\}}f_{n}(s,\rho
_{m}(Y_{s}^{n,m}),Z_{s}^{n,m})ds+K_{T}^{n,m}-K_{t}^{n,m}-%
\int_{t}^{T}Z_{s}^{n,m}dB_{s}, \\
Y_{t}^{n,m} &\geq &L_{t},\int_{0}^{T}(Y_{t}^{n,m}-L_{t})dK_{t}^{n,m}=0.
\end{eqnarray*}
Moreover, we have
\[
\left| 1_{\{s\leq \tau _{m}\}}f_{n}(s,\rho _{m}(y),z)\right| \leq m+\varphi
(m)+\beta \left| z\right| ,
\]
and $\left| \xi ^{n,m}\right| \leq m$. From Dini's theorem, we know that $%
1_{\{s\leq \tau _{m}\}}f_{n}(s,\rho _{m}(y),z)$ converge increasingly to$%
1_{\{s\leq \tau _{m}\}}f(s,\rho _{m}(y),z)$ uniformly on compact set of $%
\mathbb{R\times R}^{d}$, because $f_{n}$ are continuous and
$f_{n}$ converge increasingly to $f$. And $\xi ^{n,m}$ converge
increasingly to $\xi ^{m}$ a.s., where $\xi ^{m}=\sup_{n\geq \beta
}\xi ^{n,m}$.

As in \cite{M97}, we can prove that $Y^{n,m}$ converges increasingly to $%
Y^{m}$ in $\mathbf{S}^{2}(0,T)$, and $Z^{n,m}\rightarrow Z^{m}$ in $\mathbf{H%
}_{d}^{2}(0,T)$, $K^{n,m}\searrow K^{m}$ uniformly on $[0,T]$. Moreover, $%
(Y^{m},Z^{m},K^{m})$ solves the following RBSDE
\begin{eqnarray*}
Y_{t}^{m} &=&\xi ^{m}+\int_{t}^{T}1_{\{s\leq \tau _{m}\}}f(s,\rho
_{m}(Y_{s}^{m}),Z_{s}^{m})ds+K_{T}^{m}-K_{t}^{m}-\int_{t}^{T}Z_{s}^{m}dB_{s},
\\
Y_{t}^{m} &\geq &L_{t},\int_{0}^{T}(Y_{t}^{m}-L_{t})dK_{t}^{m}=0,
\end{eqnarray*}
where $\xi ^{m}=\sup_{n\geq \beta }Y_{\tau _{m}}^{n,m}$. Notice that $\left|
Y_{t}^{m}\right| \leq m$, so we have
\[
Y_{t}^{m}=\xi ^{m}+\int_{t}^{T}1_{\{s\leq \tau
_{m}\}}f(s,Y_{s}^{m},Z_{s}^{m})ds+K_{T}^{m}-K_{t}^{m}-%
\int_{t}^{T}Z_{s}^{m}dB_{s}.
\]
From the definition of $\{\tau _{m}\}$, it is easy to check that $\tau
_{m}\leq \tau _{m+1}$, with the definition of $Y^{m}$, $Z^{m}$, $K^{m}$ and $%
Y$, we get
\[
Y_{t\wedge \tau _{m}}=Y_{t\wedge \tau
_{m}}^{m+1}=Y_{t}^{m},Z_{t}^{m+1}1_{\{t\leq \tau
_{m}\}}=Z_{t}^{m},K_{t\wedge \tau _{m}}^{m+1}=K_{t}^{m}.
\]
We define
\[
Z_{t}:=Z_{t}^{1}1_{\{t\leq \tau _{1}\}}+\sum_{m\geq 2}Z_{t}^{m}1_{(\tau
_{m-1},\tau _{m}]}(t),\;\;K_{t\wedge \tau _{m}}:=K_{t}^{m}.
\]
Processes $(Y^{m})$ are continuous, and $P$-a.s. $\tau _{m}=T$,
for $m$ large enough, so $Y$ is continuous on $[0,T]$. It follows
that $K$ is also continuous on $[0,T]$. Furthermore, we have for
$m\in \mathbb{N}$,
\begin{equation}
Y_{t\wedge \tau _{m}}=Y_{\tau _{m}}+\int_{t\wedge \tau _{m}}^{\tau
_{m}}f(s,Y_{s},Z_{s})ds+K_{\tau _{m}}-K_{t\wedge \tau _{m}}-\int_{t\wedge
\tau _{m}}^{\tau _{m}}Z_{s}dB_{s}.  \label{st-bsde}
\end{equation}
Finally, we have
\begin{eqnarray*}
P(\int_{0}^{T}\left| Z_{s}\right| ^{2}ds &=&\infty )=P(\int_{0}^{T}\left|
Z_{s}\right| ^{2}ds=\infty ,\tau _{m}=T)+P(\int_{0}^{T}\left| Z_{s}\right|
^{2}ds=\infty ,\tau _{m}<T) \\
&\leq &P(\int_{0}^{T}\left| Z_{s}\right| ^{2}ds=\infty )+P(\tau _{m}<T),
\end{eqnarray*}
in the same way,
\[
P(\left| K_{T}\right| ^{2}=\infty )\leq P(\left| K_{T}\right| ^{2}=\infty
)+P(\tau _{m}<T).
\]
Since $\tau _{m}\nearrow T$, $P$-a.s., we know that $\int_{0}^{T}\left|
Z_{s}\right| ^{2}ds<\infty $ and $\left| K_{T}\right| ^{2}<\infty $, $P$%
-a.s. Let $m\rightarrow \infty $ in (\ref{st-bsde}), we get $(Y,Z,K)$
verifies the equation.

\textbf{Step 4}. We want to prove that the triple $(Y,Z,K)$ is a solution of
RBSDE$(\xi ,f,L)$.

First, we consider the integrability of $(Y,Z,K)$. By (\ref{est-y2}), we
know for $0\leq t\leq T$,
\begin{equation}
\left| Y_{t}\right| \leq M_{t}.  \label{est-yn}
\end{equation}
It follows immediately that
\[
E[\sup_{0\leq t\leq T}\left| Y_{t}\right| ^{2}]\leq C_{\beta }E[\left| \xi
\right| ^{2}+\int_{0}^{T}g_{s}^{2}ds+\varphi ^{2}(b)+\varphi ^{2}(2T)+1].
\]
where $C_{\beta }$ is a constant only depends on $\beta $, $T$ and $b$. For $%
K$, notice that $K^{n,m}\searrow K^{m}$, then for each $m\in \mathbb{N}$, $%
0\leq t\leq T$, we know $0\leq K_{t}^{m}\leq K_{t}^{1,m}$. Obviously, the
coefficient $1_{\{s\leq \tau _{m}\}}f_{n}(s,\rho _{m}(y),z)$ satisfies
\textbf{assumption 6'}, and Lipschitz in $z$, by Lemma \ref{est-res},
\[
E[(K_{T}^{1,m})^{2}]\leq C_{\beta }E[\left| \xi ^{1,m}\right|
^{2}+\int_{0}^{T}g_{s}^{2}ds+\varphi ^{2}(b)+\varphi ^{2}(2T)+1],
\]
where $\xi ^{1,m}=Y_{\tau _{m}}^{1}$. With (\ref{est-y2}), we have
\[
E[(K_{T}^{1,m})^{2}]\leq 2C_{\beta }E[\left| \xi \right|
^{2}+\int_{0}^{T}g_{s}^{2}ds+\varphi ^{2}(b)+\varphi ^{2}(2T)+1],
\]
which follows that for each $m\in \mathbb{N}$,
\[
E[(K_{T}^{m})^{2}]\leq 2C_{\beta }E[\left| \xi \right|
^{2}+\int_{0}^{T}g_{s}^{2}ds+\varphi ^{2}(b)+\varphi ^{2}(2T)+1],
\]
and so does for $K$, i.e. we get $E[(K_{T})^{2}]<\infty $.

In order to estimate $Z$, we apply It\^{o}'s formula to $\left| Y_{t}\right|
^{2}$ on the interval $[0,T]$, then
\begin{eqnarray*}
&&\left| Y_{0}\right| ^{2}+\frac{1}{2}E\int_{0}^{T}\left| Z_{s}\right| ^{2}ds
\\
&\leq &E\left| \xi \right| ^{2}+E\int_{0}^{T}g_{s}^{2}ds+(1+2\beta
^{2})E\int_{0}^{T}\left| Y_{s}\right| ^{2}ds+E[\sup_{0\leq t\leq T}\left|
Y_{t}\right| ^{2}]+E[(K_{T})^{2}].
\end{eqnarray*}
Thanks to the estimates for $Y$ and $K$, there exists a constant $C$ only
depends on $\beta $, $T$ and $b$, such that
\[
E\int_{0}^{T}\left| Z_{s}\right| ^{2}ds\leq CE[\left| \xi \right|
^{2}+\int_{0}^{T}g_{s}^{2}ds+\varphi ^{2}(b)+\varphi ^{2}(2T)+1].
\]
The last is to check the integral condition. Recall that $%
\int_{0}^{T}(Y_{t}^{m}-L_{t})dK_{t}^{m}=0$, then we have
\[
\int_{0}^{\tau _{m}}(Y_{t}-L_{t})dK_{t}=0,\mbox{ a.s.}
\]
Since $P$-a.s. $\tau _{m}=T$, for $m$ large enough, so
\[
\int_{0}^{T}(Y_{t}-L_{t})dK_{t}=0,\mbox{ a.s.}
\]
i.e. $(Y,Z,K)$ is a solution of RBSDE$(\xi ,f,L)$ in $\mathbf{S}%
^{2}(0,T)\times \mathbf{H}_{d}^{2}(0,T)\times \mathbf{A}^{2}(0,T)$. $\square
$

\section{Appendix: Comparison theorems}

We first generalize the comparison theorem of RBSDE with superlinear
quadratic coefficient, (in view to proposition 3.2 in \cite{KLQT}), to
compare the increasing processes. Assume that \textbf{Assumption 1} and
\textbf{3} hold, and that the coefficient $f$ satisfies:

\textbf{Assumption 7.} For all $(t,\omega )$, $f(t,\omega ,\cdot ,\cdot )$
is continuous and there exists a function $l$ strictly positive such that
\[
f(t,y,z)\leq l(y)+A\left| z\right| ^{2},\mbox{ with }\int_{0}^{\infty }\frac{%
dx}{l(x)}=+\infty .
\]

\begin{proposition}
\label{comp1}Suppose that $\xi ^{i}$ are $\mathcal{F}_{T}$-adapted and
bounded, $f^{i}(s,y,z)$, $i=1,2$ satisfy the condition \textbf{Assumption 7}
and $L$ satisfies\textbf{\ Assumption 3}. The two triples $%
(Y^{1},Z^{1},K^{1})$, $(Y^{2},Z^{2},K^{2})$ are respectively solutions of
the RBSDE$(\xi ^{1},f^{1},L)$ and RBSDE$(\xi ^{2},f^{2},L)$. If we have $%
\forall (t,y,z)\in [0,T]\times \mathbb{R\times R}^{d},
$%
\[
\xi ^{1}\leq \xi ^{2},f^{1}(t,y,z)\leq f^{2}(t,y,z),
\]
then $Y_{t}^{1}\leq Y_{t}^{2}$, $K_{t}^{1}\geq K_{t}^{2}$ and $%
dK_{t}^{1}\geq dK_{t}^{2}$, for $t\in [0,T]$.
\end{proposition}

\proof%
From the demonstration of theorem 1 in \cite{KLQT}, we know that for $i=1,2$%
, $(Y^{i},Z^{i},K^{i})$ is the solution of RBSDE$(\xi ^{i},f^{i},L)$ if and
only if $(\theta ^{i},J^{i},\Lambda ^{i})$ is the solution of RBSDE$(\eta
^{i},F^{i},\overline{L})$ where
\begin{equation}
(\theta ^{i},J^{i},\Lambda ^{i})=(\exp (2AY^{i}),2AZ^{i}\theta
^{i},2\int_{0}^{\cdot }A\exp (2AY_{s})dK_{s}^{i})  \label{trans1}
\end{equation}
and
\begin{eqnarray*}
\eta ^{i} &=&\exp (2A\xi ^{i}),\overline{L}_{t}=\exp (2AL_{t}) \\
F^{i}(t,x,\lambda ) &=&2Ax[f^{i}(s,\frac{\log x}{2A},\frac{\lambda }{2Ax})-%
\frac{\left| \lambda \right| ^{2}}{4Ax^{2}}].
\end{eqnarray*}
Then we use the approximation to construct a solution. For $p\in \mathbb{N}$%
, we consider the RBSDE$(\eta ^{i},\widetilde{F}_{p}^{i},\overline{L}_{t})$,
where
\[
\widetilde{F}_{p}^{i}(s,x,\lambda )=g(\rho (\theta ))(1-\kappa _{p}(\lambda
))+\kappa _{p}(\lambda )F^{i}(s,\rho (\theta ),\lambda ).
\]
Here $g(x)=2Axl(\frac{\log x}{2A})$, $\kappa _{p}(\lambda )$ and $\rho (x)$
are smooth functions such that $\kappa _{p}(\lambda )=1$ if $\left| \lambda
\right| \leq p$, $\kappa _{p}(\lambda )=0$ if $\left| \lambda \right| \geq
p+1$, and $\rho (x)=x$ if $x\in [r,R]$, $\rho (x)=\frac{r}{2}$ if $x\in (0,%
\frac{r}{2})$, $\rho (x)=R$ if $x\in (2R,+\infty )$, where $r$ and $R$ are
two constants. Since $\widetilde{F}_{p}^{i}$ are bounded and continuous
function of $(\theta ,\lambda )$, the RBSDE$(\eta ^{i},\widetilde{F}_{p}^{i},%
\overline{L}_{t})$ admits a bounded maximal solution $(\theta
^{i,p},J^{i,p},\Lambda ^{i,p})$, with $\underline{m}\leq $ $\theta
_{t}^{i,p}\leq V_{0}$. Here $\underline{m}$ and $V_{0}$ are constants given
in Theorem 2 in \cite{KLQT}.

We know that $\widetilde{F}_{p}^{i}\downarrow \widetilde{F}^{i}$, as $%
p\rightarrow \infty $, where $\widetilde{F}^{i}=F(s,\rho (\theta ),\lambda )$%
. Thanks to the proof of theorem 1 in \cite{KLQT}, it follows that $\theta
_{t}^{i,p}\downarrow \widetilde{\theta }_{t}^{i}$, $J_{t}^{i,p}\uparrow
\widetilde{J}_{t}^{i}$, $0\leq t\leq T$, and $\Lambda ^{i,p}\rightarrow
\widetilde{\Lambda }^{i}$ in $\mathbf{H}_{d}^{2}(0,T)$ and $(\widetilde{%
\theta }^{i},\widetilde{J}^{i},\widetilde{\Lambda }^{i})$ is a solution of
the RBSDE$(\eta ^{i},\widetilde{F}^{i},\overline{L}_{t})$. In addition, $%
\underline{m}\leq $ $\widetilde{\theta }_{t}^{i}\leq V_{0}$. So if we choose
$0<r<\underline{m}$ and $V_{0}<R$, then $\widetilde{F}^{i}=F^{i}$. It
follows that $(\widetilde{\theta }^{i},\widetilde{J}^{i},\widetilde{\Lambda }%
^{i})$ satisfies the RBSDE$(\eta ^{i},F^{i},\overline{L})$, i.e. $(%
\widetilde{\theta }^{i},\widetilde{J}^{i},\widetilde{\Lambda }^{i})=$ $%
(\theta ^{i},J^{i},\Lambda ^{i})$.

Since $f^{1}(t,y,z)\leq f^{2}(t,y,z)$, for $(t,x,\lambda )\in
[0,T]\times \mathbb{R}_{^{_{+}}}\mathbb{\times R}^{d}$, we have
$F^{1}(t,x,\lambda )\leq
F^{2}(t,x,\lambda )$. Then for $p\in \mathbb{N}$, $\widetilde{F}%
_{p}^{1}(s,x,\lambda )\leq \widetilde{F}_{p}^{2}(s,x,\lambda )$. Notice that
$\widetilde{F}_{p}^{i}$ is bounded and continuous in $(\theta ,\lambda )$
and $\theta _{t}^{i,p}>0$, by Lemma 2.1 in \cite{KLQT}, it follows that $%
\theta _{t}^{1,p}\leq \theta _{t}^{2,p}$, $J_{t}^{1,p}\geq J_{t}^{2,p}$, $%
dJ_{t}^{1,p}\geq dJ_{t}^{2,p}$, $0\leq t\leq T$. And it follows that for $%
0\leq s\leq t\leq T$, $J_{t}^{1,p}-J_{s}^{1,p}\geq J_{t}^{2,p}-J_{s}^{2,p}$.
Let $p\rightarrow \infty $, thanks to the convergence results, we get that
\[
\theta _{t}^{1}\leq \theta _{t}^{2},J_{t}^{1}\geq
J_{t}^{2},J_{t}^{1}-J_{s}^{1}\geq J_{t}^{2}-J_{s}^{2},
\]
which implies $dJ_{t}^{1}\geq dJ_{t}^{2}.$ From (\ref{trans}), we know that
\[
Y_{t}^{i}=\frac{\log (\theta _{t}^{i})}{2A},Z_{t}^{i}=\frac{\Lambda ^{i}}{%
2A\theta ^{i}},dK_{t}^{i}=\frac{dJ_{t}^{i}}{2A\theta _{t}^{i}},
\]
so $Y_{t}^{1}\leq Y_{t}^{2}$ and $dK_{t}^{1}\geq dK_{t}^{2}$, which implies
that $K_{t}^{1}\geq K_{t}^{2}$ in view of $K_{0}^{1}=K_{0}^{2}=0$. $\square $

From this result, we prove the following comparison theorem when the
coefficient $f$ satisfies monotonicity and general increasing condition in $%
y $, and quadratic increasing in $z$.

\begin{proposition}
\label{com-zq}Suppose that $\xi ^{i}$ and $f^{i}(s,y,z)$, $i=1,2$ satisfy
the condition \textbf{Asssumption 1} and \textbf{2}, $L$ satisfies \textbf{%
Assumption 3.} The two triples $(Y^{1},Z^{1},K^{1})$, $(Y^{2},Z^{2},K^{2})$
are respectively the solutions of the RBSDE$(\xi ^{1},f^{1},L)$ and RBSDE$%
(\xi ^{2},f^{2},L)$. If we have $\forall (t,y,z)\in [0,T]\times \mathbb{%
R\times R}^{d}$,
\[
f^{1}(t,y,z) \leq f^{2}(t,y,z), \quad \xi ^{1} \leq \xi ^{2},
\]
then $Y_{t}^{1}\leq Y_{t}^{2}$, $K_{t}^{1}\geq K_{t}^{2}$ and $%
dK_{t}^{1}\geq dK_{t}^{2}$, for $t\in [0,T]$.
\end{proposition}

\proof%
First with changement of $(Y,Z,K)$,
\[
(Y^{b},Z^{b},K^{b})=(Y-b,Z,K),
\]
we work with $L^{b}\leq 0$. Since this transformation doesn't
change the monotonicity, then in the following, we assume that the
barrier $L$ is a negative bounded process. As in the proof of
theorem \ref{gmq}, for $C\in \mathbb{R}_{+}$, let
$g^{C}:\mathbb{R}\rightarrow [0,1]$ continuous which
satisfies (\ref{gfun}). Set $f_{i}^{C}(t,y,z)=g^{C}(x)f^{i}(t,y,z)$, $i=1,2$%
, which satisfies \textbf{Assumption 7}, with $l^{i}(y)=\varphi ^{i}(\left|
2C\right| )$. We consider solutions $(Y^{i,C},Z^{i,C},K^{i,C})$ of the RBSDE$%
(\xi ^{i},f_{i}^{C},L^{b})$ respectively. Using proposition \ref{comp1},
since
\[
f_{1}^{C}(t,y,z)\leq f_{2}^{C}(t,y,z),\mbox{ }\xi ^{1}\leq \xi ^{2},
\]
we get for $t\in [0,T]$,
\[
Y_{t}^{1,C}\leq Y_{t}^{2,C},dK_{t}^{1,C}\geq dK_{t}^{2,C}.
\]
Then by the bounded property of $Y^{i}$, we choose $C$ big enough like in
the proof of theorem \ref{gmq}, which follows immediately
\[
Y_{t}^{1}\leq Y_{t}^{2},dK_{t}^{1}\geq dK_{t}^{2},\forall t\in [0,T].
\]
$\square $

\begin{proposition}
\label{com-zl1}Suppose that $\xi ^{i}\in \mathbf{L}^{2}\mathbf{(}\mathcal{F}%
_{T})$, $f^{i}(s,y,z)$ satisfy the condition \textbf{Assumption 6}, and $%
L^{i}$ satisfies \textbf{Assumption 3, }$i=1,2$. The two triples $%
(Y^{1},Z^{1},K^{1})$, $(Y^{2},Z^{2},K^{2})$ are respectively solutions of
the RBSDE$(\xi ^{1},f^{1},L)$ and RBSDE$(\xi ^{2},f^{2},L)$. If we have for $%
\forall (t,y,z)\in [0,T]\times \mathbb{R\times R}^{d},$%
\[
\mbox{ }\xi ^{1}\leq \xi ^{2},\,\,\,f^{1}(t,y,z)\leq
f^{2}(t,y,z),\,L_{t}^{1}\leq L_{t}^{2},
\]
then $Y_{t}^{1}\leq Y_{t}^{2}$, for $t\in [0,T]$.
\end{proposition}

The result comes from the comparison theorem in \cite{LMX} and the
approximation in the proof of theorem \ref{exist-zl}.


\begin{thebibliography}{99}
\bibitem{BDHPS}  Briand, Ph., Delyon, B., Hu, Y., Pardoux, E. and Stoica, L.
(2003) $L_{p}$ solutions of BSDEs, \textit{Stochastic Process. Appl.} 108,
109-129.

\bibitem{LS04}  Ph. Briand, J.P. Lepeltier and J. San Mart\'{i}n, (2004)
BSDE's with continuous, monotonicity, and non-Lipschitz in $z$ coefficient.
Submitted

\bibitem{E79}  N. El Karoui, (1979) Les aspects probabilistes du contr\^{o}%
le stochastique. Ecole d'\'{e}t\'{e} de Saint-Flour, \textit{Lecture Notes
in Math. }\textbf{876. (}Springer, Berlin), 73-238.

\bibitem{EKPPQ}  N.~El Karoui, C.~Kapoudjian, E.~Pardoux, S.~Peng and M.C.
Quenez, (1997) Reflected Solutions of Backward SDE and Related Obstacle
Problems for PDEs, \textit{Ann. Probab.} \textbf{25}, no 2, 702--737.

\bibitem{EPQ}  N.~El Karoui, S.~Peng and M.C. Quenez, (1997). Backward
stochastic differential equations in Finance. \textit{Math. Finance, }%
\textbf{7, }1-71\textbf{.}

\bibitem{K00}  M. Kobylanski, (2000) Backward stochastic differential
equations and partial differential equations with quadratic growth, \textit{%
Ann. Proba. }\textbf{28}, 558-602.\textbf{\ }

\bibitem{KLQT}  M. Kobylanski, J.P. Lepeltier, M.C. Quenez and S.Torres,
(2002) Reflected BSDE with superlinear quadratic coefficient. \textit{%
Probability and Mathematical Statistics. }Vol. \textbf{22}, 51-83.

\bibitem{LMX}  J.P. Lepeltier, A. Matoussi and M. Xu, (2004) Reflected BSDEs
under monotonicity and general increasing growth conditiond. \textit{%
Advanced in applied probability}, to appear in March, 2005.

\bibitem{LS98}  J.P. Lepeltier and J. San Mart\'{i}n. (1998) Existence for
BSDE with superlinear-quadratic coefficient, \textit{Stochastic and
Stochastic reports,} \textbf{63}, 227-240.

\bibitem{M97}  A. Matoussi, (1997) Reflected solutions of backward
stochastic differential equations with continuous coefficient, \textit{%
Statistic \& Probality Letters} \textbf{34}, 347-354.

\bibitem{P99}  Pardoux, E., 1999. BSDE's, weak convergence and
homogenization of semilinear PDE's in Nonlinear analysis, Differential
Equations and Control, F. H. Clarke \& R. J. Stern Eds , pp. 503-549, Kluwer
Acad. Pub.

\bibitem{PP90}  E. Pardoux and S. Peng, (1990) Adapted solutions of Backward
Stochastic Differential Equations. \textit{Systems Control Lett. }\textbf{14,%
} 51-61.

\bibitem{RY}  D. Revuz and M. Yor, (1991) Continuous artingales and Brownian
motion (Springer, Berlin).
\end{thebibliography}
\end{document}